\documentclass[11pt]{amsart}
\usepackage{amsmath}
\usepackage{amssymb}
\usepackage{stackrel,amssymb}
\usepackage{amsfonts}
\usepackage{amsthm}
\usepackage{slashed}
\usepackage{tikz-cd}
\usepackage{mathrsfs}
\usepackage[all]{xy}
\usepackage{mathtools}

\usepackage{mathabx,epsfig}

\usepackage{xcolor}

\usepackage{amssymb,amsfonts}
\usepackage{enumerate}
\usepackage{enumitem}
\usepackage{mathrsfs}
\usepackage{stmaryrd}
\usepackage{hyperref}
\usepackage{cleveref}

\usepackage{relsize}
\usepackage[bbgreekl]{mathbbol}
\usepackage{amsfonts}
\DeclareSymbolFontAlphabet{\mathbb}{AMSb} 
\DeclareSymbolFontAlphabet{\mathbbl}{bbold}

\begin{document}

\newcommand\restr[2]{{
  \left.\kern-\nulldelimiterspace 
  #1 
  \vphantom{\big|} 
  \right|_{#2} 
  }}

\newtheorem{thm}{Theorem}[subsection]
\newtheorem{cor}[thm]{Corollary}
\newtheorem{prop}[thm]{Proposition}
\newtheorem{lem}[thm]{Lemma}
\newtheorem{conj}[thm]{Conjecture}
\newtheorem{quest}[thm]{Question}

\theoremstyle{definition}
\newtheorem{defn}[thm]{Definition}
\newtheorem{defns}[thm]{Definitions}
\newtheorem{con}[thm]{Construction}
\newtheorem{exmp}[thm]{Example}
\newtheorem{exmps}[thm]{Examples}
\newtheorem{notn}[thm]{Notation}
\newtheorem{notns}[thm]{Notations}
\newtheorem{addm}[thm]{Addendum}
\newtheorem{exer}[thm]{Exercise}

\theoremstyle{remark}
\newtheorem{rem}[thm]{Remark}
\newtheorem{rems}[thm]{Remarks}
\newtheorem{warn}[thm]{Warning}
\newtheorem{sch}[thm]{Scholium}

\makeatletter
\renewcommand{\@seccntformat}[1]{%
  \ifcsname prefix@#1\endcsname
    \csname prefix@#1\endcsname
  \else
    \csname the#1\endcsname\quad
  \fi}
\makeatother

\title{Character Formulas from Matrix Factorizations}         
\author{Kiran Luecke}        
\date{\today}          
\maketitle

\begin{abstract}
In this paper I present a new and unified method of proving character formulas for discrete series representations of connected Lie groups by applying a Chern character-type construction to the matrix factorizations of [FT] and [FHT3]. In the case of a compact group I recover the Kirillov formula, thereby exhibiting the work of [FT] as a categorification of the Kirillov correspondence. In the case of a real semisimple group I recover the Rossman character formula with only a minimal amount of analysis. The appeal of this method is that it relies almost entirely on highest-weight theory, which is a far more ubiquitous phenomenon than the varied techniques that were previously used to prove such formulas.

\end{abstract}


\tableofcontents

\section*{Introduction}

Let $G$ be a compact Lie group. The famous Kirillov correspondence establishes a close relationship between the representation theory of $G$ and certain coadjoint orbits in the dual of the Lie algebra. On the other hand, a celebrated theorem of Freed, Hopkins, and Teleman [FHT3] establishes a close relationship between the positive energy representation theory of the loop group $LG$ and the twisted, conjugation-equivariant $K$-theory of $G$. The statement is cleanest when $G$ is simple and simply-connected with dual Coxeter number $h$. Then $H^3_{G}G\simeq\mathbb{Z}$ canonically, so an integer $\tau$ represents both a level for positive energy representations as well as a twist for equivariant $K$-theory. Let $PER^\tau(LG)$ denote the category of positive energy representations of level $\tau$. Then there is an isomorphism
$$^{\tau+h} K^{\text{dim}G}_GG\simeq  K_0PER^\tau(LG).$$
A calculation of the left side shows that the $K$-classes are supported at a discrete set of conjugacy classes ([FHT3] Proposition 6.12). Since conjugacy classes in $G$ correspond to coajdoint orbits of $LG$, one suspects that the Kirillov correspondence is at play. The isomorphism, although established by calculating the left side, is explained by a map from right to left which constructs out of a positive energy representation on the right a twisted family of Fredholm operators representing a class on the left ([FHT3] Section 4). This construction crucially uses a certain family of Dirac operators. In [FT] Freed and Teleman re-interpret this Fredholm family as a matrix factorization, allowing them to categorify the isomorphism of abelian groups in the previous display. Moreover, the construction applies (in a much simpler form) to the coadjoint representation, providing an explicit map
$$K^{\text{dim}\mathfrak{g}}_{G}(\mathfrak{g}^{*})\xleftarrow{\sim} K^0_G(\text{pt})\simeq R(G)$$
presenting the Thom isomorphism. Combined with the categorical perspective of [FT] and the Kirillov philosophy, one has arrived at a categorification of the Kirillov correspondence, and indeed, the main theorems in [FHT3] and [FT] should be interpreted that way. 

In this paper I confirm this philosophy by extracting the Kirillov character formula from the categorical correspondence of [FT]. The method of extraction is a Chern character-type construction that applies to matrix factorizations quite generally. In particular Frenkel's character formula for loop groups likely follows from a proper application of this method to the categorification of the first displayed isomorphism above\footnote{However, the additional amount of analysis required for the case of $LG$ is probably proportional to the difference in dimension between $LG$ and $G$.} as well as character formulas for more general Kac-Moody groups [Ki]. However, in this paper I generalize in a different direction, to finite-dimensional but non-compact groups by proving the Rossman character formula for discrete series representations of connected real semisimple Lie groups. In the remainder of this section I recall the character formulas and give a brief outline of the proofs that make up the rest of the paper.


Let $G$ be a connected real semisimple Lie group. Recall that the \textit{discrete series} representations are the family\footnote{As the name suggests, this family is indexed by discrete parameters.} of (isomorphism classes of) unitary representations of $G$ defined by the property that all matrix coefficients are in $L^{2}(G)$. By a well-known result of Harish-Chandra the discrete series is empty unless the rank of $G$ is equal to the rank of its maximal compact subgroup. Fix an irreducible discrete series representation ($\mathcal{H}_{-\Lambda},\pi)$. Let $U$ be an open neighborhood of $0\in\mathfrak{g}$ over which $\text{exp}: U\rightarrow G$ is invertible and analytic. For a compactly supported smooth function $\varphi\in C^{\infty}_{c}(U)$ define the smeared operator $\pi(\varphi):=\int_{\mathfrak{g}}dX\ \varphi(X)\pi(e^{X})$. Define the equivariant $\hat{\mathcal{A}}$-genus at an element $X\in\mathfrak{g}$ as follows: let $P$ be a system of positive roots. Then
$$\hat{\mathcal{A}}(X):=\prod_{\alpha\in P}\frac{e^{\langle\alpha,X\rangle/2}-e^{-\langle\alpha,X\rangle/2}}{\langle\alpha,X\rangle}.$$
The operator $\pi(\varphi)$ is trace class and satisfies the following formula.

\begin{thm}[Rossman] 
$$\text{Tr}_{\mathcal{H}_{-\Lambda}}\Big(\int_{\mathfrak{g}}dX\  \varphi(X)\pi(e^{X})\Big)=\int_{\mathfrak{g}^{*}}\delta_{O_{-\Lambda-\tilde{\rho}}}(\xi)\ \Big(\int_{\mathfrak{g}} e^{i\langle\xi,X\rangle}\varphi(X)\hat{\mathcal{A}}(X)dX\Big)d\xi$$
\end{thm}

If $G$ is compact then every irreducible representation is in the discrete series and individual group elements act by trace class operators, so the above formula simplifies to the well-known Kirillov formula.

\begin{thm}[Kirillov] Let $G$ be a compact connected Lie group and $V_{-\lambda}$ the irreducible representation of lowest weight $-\lambda$. For $X\in\mathfrak{g}$ close to 0
$$\text{Tr}_{V_{-\lambda}}(\pi(e^{X}))=\hat{\mathcal{A}}(X)\int_{\mathfrak{g}^{*}}\delta_{O_{-\lambda-\rho}}(\xi)e^{i\langle\xi,X\rangle}d\xi.$$
\end{thm}

These formulas will be proved as follows. From the data of the matrix factorizations of [FHT3] and [FT] one constructs a super-bundle over $\mathfrak{g}^*$ together with a one-parameter family of super-connections. That leads (via the Mathai-Quillen formula) to a one-parameter family of differential forms representing the Chern character of the super-bundle. Using highest weight theory and knowledge of the spectral properties of the Dirac family defining the matrix factorization these differential forms can be rewritten in such a way that the one-parameter family interpolates between the integrands on either side of the Kirillov formula. The equality is then established by appealing to the crucial fact that the cohomology class of the representative of the Chern character is invariant under one-parameter deformations. It is worth noting that in the case of compact groups, one can prove the Kirillov formula by appealing to the index theorems in equivariant $K$-theory. However the method outlined above is more robust in the sense that it continues to apply when $G$ is not compact and $K$-theoretic tools are not available.

I would like to thank Constantin Teleman for helpful conversations and comments on earlier drafts of this paper, Mikayla Kelley for guidance regarding some integral estimates, and an anonymous referee for valuable comments and corrections.

\section{The Kirillov Formula}

\textit{Throughout this section, $G$ is a compact connected Lie group and $(V,\pi)$ is an irreducible representation of lowest weight $-\lambda$.}

Choose a $G$-invariant inner product on $\mathfrak{g}^*$ (e.g. the Killing form) and let $\mathcal{S}$ be an irreducible $\mathbb{Z}/2$-graded Cliff$^{c}(\mathfrak{g}^{*})$-module\footnote{See [FHT2 Prop 1.18]}. Consider the equivariant super-bundle $\mathcal{V}:=V\otimes\mathcal{S}\times\mathfrak{g}^{*}\rightarrow\mathfrak{g}^{*}$. Since $\mathcal{S}$ carries a projective action of $G$, this bundle is projectively $G$-equivariant. Let $\gamma(\xi)$ denote Clifford multiplication by $\xi\in\mathfrak{g}^{*}$, and let $\iota(X)$ denote contraction with the vector field defined by $X\in\mathfrak{g}$ acting in the coadjoint representation. Let $R(X)$ and $\sigma(X)$ denote the Lie algebra action on $V$ and $\mathcal{S}$, respectively. Let $\{e^{\nu}\}$ be a basis of $\mathfrak{g}^{*}$ with dual basis $\{e_{\nu}\}$ and abbreviate $R(e_{\nu})$, $\sigma(e_{\nu})$, and $\gamma(e^{\nu})$ by $R_{\nu}$, $\sigma_{\nu}$, and $\gamma^{\nu}$. Let $\slashed D_{0}$ Kostant's cubic Dirac operator (c.f. [FHT2] 1.10), the endomorphism of $V\otimes\mathcal{S}$ given by\footnote{This differs from the operator used in [FHT1-3] by a factor of $i$.}

$$\slashed D_{0}=\sum_{\nu}\big(R_{\nu}+\frac{1}{3}\sigma_{\nu}\big)\otimes\gamma^{\nu}.$$ 
Recall that a \textit{super-connection} on $\mathcal{V}$ is an odd\footnote{Note that the total grading on the two complexes above---with respect to which the requirement of oddity is made---is the sum of the grading mod 2 on differential forms and the super-grading on $\mathcal{V}^{\pm}$.} differential operator on bundle-valued forms
$$A:\Omega^{\bullet}(\mathfrak{g}^{*},\mathcal{V})^{\pm}\longrightarrow\Omega^{\bullet}(\mathfrak{g}^{*},\mathcal{V})^{\mp}.$$
For $X\in\mathfrak{g}$ let $\iota(X)$ denote contraction with the vector field on $\mathcal{V}$ defined by $X$. Then the \textit{equivariant refinement} of $A$ is an operator on equivariant forms

$$ A:(\mathbb{C}[\mathfrak{g}]\otimes\Omega^{\bullet}(\mathfrak{g}^{*},\mathcal{V})^{G})^{\pm}\longrightarrow(\mathbb{C}[\mathfrak{g}]\otimes\Omega^{\bullet}(\mathfrak{g}^{*},\mathcal{V})^{G})^{\mp}$$
namely the $\Omega^{\bullet}(\mathfrak{g}^{*},\mathcal{V})^{G})$-valued function on $\mathfrak{g}$ defined by the formula
$$A(X):=A -\iota(X).$$
Let $\mathsterling(X)$ be the Lie derivative along the vector field on $\mathcal{V}$ defined by $X$ and define the \textit{moment} of $X$ as $\vec{\mu}(X):=\mathsterling(X)-[\iota(X), A]$. The \textit{equivariant Chern character} of $\mathcal{V}$ with respect to $A$, for a test element $X\in\mathfrak{g}$, is defined following Mathai-Quillen to be the differential form

$$\text{str}(e^{A(X)^2})=\text{str}(e^{A^2+\vec{\mu}(X)}).$$
\begin{rem}\label{topdeg}
Although the ``differential form" above is, strictly speaking, a formal sum of differential forms of different degree (one for each term in the Taylor series of $e^x$), a single term of homogeneous degree is singled out when integrated over a manifold of fixed dimension.
\end{rem}
To prove the Kirillov character formula I make heavy use of the following result, which is Theorem 7.7 of [BGV].

\begin{lem}\label{deformation} Under one parameter deformations of the super-connection, the equivariant Chern character differs by an equivariantly exact differential form.
\end{lem}

\begin{proof} Indeed, suppose $A_{\epsilon}(X)$ is a one parameter family of equivariant superconnections. Let $d_{G}$ denote the equivariant de Rham operator. Then, making use of the equivariant Bianchi identity $[A_{\epsilon}(X),A_{\epsilon}(X)^2]=0$,
\begin{align*}
\frac{d}{d\epsilon}\text{str}(e^{A_{\epsilon}(X)^2}) & =\text{str}\big([A_{\epsilon}(X), \frac{dA_{\epsilon}}{d\epsilon}(X)]e^{A_{\epsilon}(X)^2}\big)\\
&=\text{str}\big([A_{\epsilon}(X), \frac{dA_{\epsilon}}{d\epsilon}(X)e^{A_{\epsilon}(X)^2}]\big)\\ 
& =d_{G}\Big(\text{str}(\frac{dA_{\epsilon}}{d\epsilon}(X)e^{A_{\epsilon}(X)^2})\Big).
\end{align*}
\end{proof}

\begin{defn}\label{family}
Consider the following 2-parameter family of equivariant super-connections on $\mathcal{V}$.

$$A_{\epsilon}^{t}(X):=d +\sqrt{t}\slashed D_{0}+  i\sqrt{\epsilon}\gamma(\xi)-\iota(X).$$
\end{defn}

\begin{lem}\label{decay}
When $\epsilon\neq0$ the differential form $\text{str}(e^{A_{\epsilon}^{t}(X)^2})$ has Gaussian decay as $\xi\to\infty$.
\end{lem}
\begin{proof}
Expanding $A_{\epsilon}^{t}(X)^2$ using Definition \ref{family} shows that the leading term in $\xi$ is the quadratic term $-\epsilon\gamma(\xi)^2=-\epsilon|\xi|^2$, all other terms are either linear or constant in $\xi$.
\end{proof}
  
\begin{rem}
Although the definition of $\mu(X)$ depends on a choice of connection, for any connection in the family  $A_{\epsilon}^{t}(X)$ defined above, the moment $\mu(X)$ is the same. That's because $\iota(X)$ commutes with $\slashed D_{0}$ and $\gamma(\epsilon)$ (viewing a differential form as a function from some power of the tangent space to the fiber of the bundle $\mathcal{V}$, contraction is a precomposition and $\slashed D_{0}$ and $\gamma(\epsilon)$ are post-composition) and $A_{\epsilon}^{t}(X)$ differs from the trivial connection by linear combinations of those.
\end{rem}

Let $t=0$. The equivariant Chern character of the remaining $\epsilon$-family is 

$$\text{str}(e^{A_{\epsilon}^{0}(X)^2})=\text{str}(e^{-\epsilon|\xi|^2+i\sqrt{\epsilon}d\gamma(\xi)+ \vec{\mu}(X)}).$$
Lemma \ref{deformation}, Stokes theorem, and the decay established in Lemma \ref{decay} imply that for any $\epsilon_{1},\epsilon_{2}>0$,

\begin{align*}
\int_{\mathfrak{g}^{*}}\text{str}(e^{A^{0}_{\epsilon_{1}}(X)^2}-e^{A^{0}_{\epsilon_{2}}(X)^2})d\xi &=\lim_{r\to\infty}\int_{|\xi |=r}d\xi\int_{\epsilon_{1}}^{\epsilon_{2}}\text{str}(\frac{i}{2\sqrt{\epsilon}}\gamma(\xi)e^{A^{0}_{\epsilon}(X)^2})d\epsilon\\
&=0.
\end{align*}
As $\epsilon\to\infty$, the differential form $\text{str}(e^{A_{\epsilon}^{0}(X)^2})$ localizes to the kernel of $\gamma(\xi)^{2}=-|\xi|^2$, which is the fiber of $\mathcal{V}$ over $\xi=0$. But there the moment $\vec{\mu}(X)$ is simply the Lie algebra action of $X$ in the representation $V\otimes\mathcal{S}$. Recall that the \textit{chirality operator} of the Clifford algebra Cliff$^c(\mathfrak{g}^*)$ is defined to be the product $\Gamma=\gamma(e_1)...\gamma(e_n)$ for any choice of orthonormal basis $e_i$ of $\mathfrak{g}^*$. The supertrace $\text{str}_\mathcal{S}(\Gamma)$ is nonzero (c.f. [BGV] Lemma 3.17). Moreover, the degree $n$ part of the supertrace over $\mathcal{S}$ of $e^{id\gamma(\xi)+ \vec{\mu}(X)}$ is $\text{str}_\mathcal{S}(\Gamma)\hat{\mathcal{A}}(X)^{-1}$ (c.f. [BGV] Propositions 3.16 and 3.21). In conclusion, 

\begin{align*}
\lim_{t=0,\epsilon\to\infty}\text{str}(e^{A_{\epsilon}^{t}(X)^2})&=\lim_{t=0,\epsilon\to\infty}\text{str}_\mathcal{S}(\Gamma)\epsilon^{\frac{n}{2}}e^{-\epsilon|\xi|^2}d\xi\hat{\mathcal{A}}^{-1}(X)\text{Tr}_{V}(e^{X})\\
&=\text{str}_\mathcal{S}(\Gamma)\pi^{\frac{n}{2}}\delta_{0}(\xi)d\xi\hat{\mathcal{A}}^{-1}(X)\text{Tr}_{V}(e^{X}).
\end{align*}
Now let $t=\epsilon$. Write $\xi=\sum_{\nu}\xi_{\nu}e^{\nu}$ and define $T(\xi):=\sum_{\nu}\xi_{\nu}(R_{\nu}+\sigma_{\nu})$. Then the Chern character becomes

$$\text{str}(e^{A_{\epsilon}^{\epsilon}(X)^2})=\text{str}(e^{ \epsilon\slashed D_{0}^2+2i\epsilon T(\xi) -\epsilon|\xi|^2+i\sqrt{\epsilon}d\gamma(\xi)+\vec{\mu}(X)}).$$
By the same argument as above, as $\epsilon\to\infty$ this differential form localizes to the kernel of $(\slashed D_{0}+i\gamma(\xi))^2=\slashed D_{0}^2+2i T(\xi) -|\xi|^2$. Let $\mathcal{S}_{N}$ denote the spinor bundle of the normal bundle $N$ of $O_{-\lambda-\rho}\hookrightarrow\mathfrak{g}^{*}$. By [FHT2] Proposition 1.19, the operator $(\slashed D_{0}+i\gamma(\xi))^2\in\text{End}(V\otimes\mathcal{S})$ is negative semi-definite and diagonalizable along a weight basis corresponding to a Cartan subalgebra containing the dual of $\xi$, has a maximum eigenvalue of $|\lambda+\rho+\xi|^2$ which occurs on the lowest weight space, which traces out a sub-bundle isomorphic to $\mathcal{L}(-\lambda-\rho)\otimes\mathcal{S}_{N}$ as $\xi$ varies (and the Cartan subalgebra defining the weigh basis along with it) along the orbit $O_{-\lambda-\rho}$ (c.f. [FHT2] Proposition 1.19). Let $\mathfrak{h}$ be the Cartan subalgebra stabilizing the element $\lambda+\rho\in O_{-\lambda-\rho}$. Let $\Lambda(V\otimes\mathcal{S})$ be the graded set of weights (with multiplicities, and with respect to $\mathfrak{h}$) of the graded representation $V\otimes\mathcal{S}$ and let $|\tau_{\pm}\rangle$ be a weight basis labelled by $\tau_{\pm}\in\Lambda(V\otimes\mathcal{S})$. By the facts mentioned above, for any $\xi\in\mathfrak{h}^*$ there are non-negative real numbers $c(\tau_\pm)$ such that the operator  $(\slashed D_{0}+i\gamma(\xi))^2$ has an eigenvalue of $|\lambda+\rho+\xi|^2-c(\tau_\pm)$ on $|\tau_{\pm}\rangle$. Note that $\mathfrak{h}$ defines a transverse slice to $O_{-\lambda-\rho}$ at $\lambda+\rho$, namely the subspace $\mathfrak{h}^*\subset \mathfrak{g}^*$. Let $e^{1},...,e^{n}$ be a basis of $\mathfrak{g}^*$ such that $e^{1},...,e^{r}$ are a basis of $\mathfrak{h}^*$. Note that the tangent space $T_{\lambda+\rho}O_{-\lambda-\rho}$ defines a subspace of $\mathfrak{g}^*$ complimentary to $\mathfrak{h}^*$. Then contraction of the top degree part of $\text{str}(e^{A_{\epsilon}^{\epsilon}(X)^2})$ (c.f. Remark \ref{topdeg}) with a basis of $T_{\lambda+\rho}O_{-\lambda-\rho}$ $\mathfrak{h}^*$ can be written as

$$e^{-\epsilon|\lambda+\rho+\xi|^2}\epsilon^{\frac{r}{2}}\text{str}_\mathcal{S}(\Gamma)d\xi^1\wedge...\wedge\ d\xi^{r}\sum_{\tau_{\pm}\in\Lambda(V\otimes\mathcal{S})}\pm e^{\epsilon c(\tau_{\pm})}\langle\tau_{\pm}|e^{\vec{\mu}(X)}|\tau_{\pm}\rangle.$$
The differential form $\text{str}(e^{A_{\epsilon}^{\epsilon}(X)^2})$ is equivariant (and so is the decomposition of $\mathfrak{g}^*$ into $\mathfrak{h}^*$ and $T_{\lambda_\rho}O_{-\lambda-\rho}$) so the formula above determines it along the entire coadjoint orbit. Therefore,

$$\lim_{t=\epsilon\to\infty}\text{str}(e^{A_{\epsilon}^{t}(X)^2})=\text{str}_\mathcal{S}(\Gamma)\pi^{\frac{r}{2}}\delta_{-\lambda-\rho}(\xi)d\xi\langle-\lambda-\rho|e^{\vec{\mu}(X)}|-\lambda-\rho\rangle.$$
As mentioned above bundle spanned by $|-\lambda-\rho\rangle$ as $\xi$ varies along $O_{-\lambda-\rho}$ is exactly $\mathcal{L}(-\lambda-\rho)$, on which the supertrace of $e^{\vec{\mu}(X)}$ is (cohomologous to) $e^{i\langle X,\xi\rangle}$ (c.f. [BGV] Proposition 8.6). Note that the moment on spinors is trivial since the normal bundle to $O_{-\lambda-\rho}$ at any point is contained in the stabilizer of that point. Therefore
$$\lim_{t=\epsilon\to\infty}\text{str}(e^{A_{\epsilon}^{t}(X)^2})=\text{str}_\mathcal{S}(\Gamma)\delta_{O_{-\lambda-\rho}}(\xi)e^{i\langle\xi,X\rangle}d\xi.$$

Choose a smooth path in the $(\epsilon,t)$-plane which tends to  $(\infty,0)$ at one end and becomes $\epsilon=t\to\infty$ at the other. Apply Lemma \ref{deformation} to this one parameter family of super-connections to complete the proof of the Kirillov formula.


 \section{The Dirac Family for the Discrete Series}
Let $G$ be a connected real semisimple Lie group with maximal compact subgroup $K$ and Cartan decomposition $\mathfrak{g}=\mathfrak{k}\oplus\mathfrak{p}$. Choose a maximal torus $T\subset K$ with Lie algebra $\mathfrak{h}$ and assume that $T$ is a Cartan subgroup of $G$ so that the discrete series is non-empty. Throughout this section $\mu\in\mathfrak{h}^{*}\subset\mathfrak{k}^{*}$ will be a regular elliptic element. Choose the system of positive roots for which $\mu$ is antidominant; it admits a decomposition $\Delta^{+}=\Delta^{+}_{c}\cup\Delta^{+}_{n}$ into compact roots (contained in $\mathfrak{k}^*$) and noncompact roots (contained in $\mathfrak{p}^*$. Write $\rho_c:=\frac{1}{2}\sum_{\Delta_c^+}\alpha$ and $\rho_n:=\frac{1}{2}\sum_{\Delta_n^+}\alpha$ and put $\tilde{\rho}:=\rho_{c}-\rho_{n}$. Denote the Killing form by $\langle\cdot,\cdot\rangle$. Let $\mathcal{S}$ be an irreducible $\mathbb{Z}/2$-graded spinor representation of Cliff$^{c}(\mathfrak{g}^{*},\langle\cdot,\cdot\rangle)$ and let $S_{-\tilde{\rho}}$ be its lowest weight line. Choose a basis $\{e_{\nu}\}$ of $\mathfrak{g}^{\mathbb{C}}$ with dual basis $\{e^{\nu}\}$ which is compatible with both the chosen root and Cartan decomposition. For a subspace $\mathfrak{s}\subset\mathfrak{g}$ write $\nu\in\mathfrak{s}$ when $e_{\nu}\in\mathfrak{s}$. As before, $\gamma^{v}$ denotes Clifford multiplication by $e^{\nu}$, $R_{\nu}$ and $\sigma_{\nu}$ denote the action of $e_{\nu}$ on $\mathcal{H}_{\Lambda}$ and $\mathcal{S}$, and $T(\mu)=\sum_{\nu}\mu_{\nu}(R_{\nu}+\sigma_{\nu})$.  Recall that discrete series representations enjoy the property of having a dense subspace which is a direct sum of irreducible representations of $K$, called \textit{$K$-types}, each occurring with at most finite multiplicity. Furthermore, there is a distinguished $K$-type which occurs with multiplicity one and is \textit{minimal}---if $-\Lambda$ denotes its lowest weight and $-\Lambda'$ denotes the lowest weight of any other $K$-type, then $|\Lambda+\tilde{\rho}|^2<|\Lambda'+\tilde{\rho}|^2$. Let $\mathcal{H}_{\Lambda}$ be a discrete series representation with minimal $K$-type $-\Lambda$ and denote by $L_{-\Lambda}$ its lowest weight line.

\begin{lem} The negative non-compact roots elements annihilate $L_{-\Lambda}$.
\end{lem}

\begin{proof} Let $v\in L_{-\Lambda}$. Take $X_{-\beta}\in\mathfrak{g}_{-\beta}$ a root element for some $\beta\in\Delta^{+}_{n}$. Then $X_{-\beta}\cdot v = 0$ or $X_{-\beta}\cdot v\in L_{-\Lambda -  \beta}$. In the latter case, by applying some (possibly empty) combination of compact roots $\Delta_{c}$ one may conclude that $X_{-\beta}\cdot v$ is contained in a $K$-type of highest weight $$-\Lambda - \beta + \sum\limits_{\alpha\in\Delta^{+}_{c}}n_{\alpha}\alpha;  (n_{\alpha}\geq0).$$
But by [Kn] Theorem 9.20 the only possible $K$-types in such a setting are of the form $$-\Lambda^{'}=-\Lambda+\sum_{\nu\in\Delta^{+}\cup\Delta^{-}_{c}}n_{\nu}\nu; (n_{\nu}\geq0).$$
Hence it must be that $X_{-\beta}\cdot v = 0$, as desired, because the above formulas disagree on the sign of the non-compact simple root, which appears an odd number of times in $\beta$ and an even number of times in any compact root. 
\end{proof}

\begin{defn}
For each $\mu=\sum\mu_{\nu}e^{\nu}\in\mathcal{E}^{*}$ define three endomorphisms of $\mathcal{H}_{\Lambda}\otimes\mathcal{S}$
$$\big(\slashed D_{\mu}\big)_{c}:= \sum_{\nu\in\mathfrak{k}}(R_{\nu}+\frac{1}{3}\sigma_{\nu} + i\mu_{\nu})\otimes\gamma^{\nu},$$
$$\big(\slashed D_{\mu}\big)_{n}:= \sum_{\nu\in\mathfrak{p}}(R_{\nu}+\frac{1}{3}\sigma_{\nu} + i\mu_{\nu})\otimes\gamma^{\nu},$$
 $$\slashed D_{\mu}:= i\big(\slashed D_{\mu}\big)_{c}+\big(\slashed D_{\mu}\big)_{n}.$$ 
\end{defn}

\begin{defn}\label{locus}
Define $\mathcal{E}^*$ to be the subset of $\mathfrak{g}^*$ consisting of elements of positive square norm, as measured by the Killing form. Note that elliptic elements are contained in $\mathcal{E}^*$.
\end{defn}  
\begin{prop} The operator $\slashed D_{\mu}^2\in\text{End}(\mathcal{H}_{\Lambda}\otimes\mathcal{S})$ is invertible for all elements $\mu\in\mathcal{E}^*$ except on the coadjoint orbit $O_{-\Lambda-\tilde{\rho}}$. For $\mu$ in this orbit, 
$\text{ker}\slashed D_{\mu}=\text{ker}\slashed D_{\mu}^2=L_{-\Lambda}\otimes S_{-\tilde{\rho}}$. 
Away from $O_{-\Lambda-\tilde{\rho}}$ the spectrum of $\slashed D_{\mu}^2$ is gapped a finite distance away from 0. There is an associated family of bounded Fredholm operators $\slashed D_{\mu}(1-\slashed D_{\mu}^2)^{-1/2}$ acting on the Hilbert space $\mathcal{H}_{\Lambda}\otimes \mathcal{S}$ which has all the spectral properties stated for $\slashed D_{\mu}^2$.
\end{prop}

\begin{proof} Let $\mathfrak{n}_{c}=\bigoplus_{\alpha\in\Delta^{+}_{c}}\mathfrak{g}^{\mathbb{C}}_{\alpha}$, $\mathfrak{n}_{n}=\bigoplus_{\alpha\in\Delta_{n}^{+}}\mathfrak{g}^{\mathbb{C}}_{\alpha}$ and write $\mathfrak{n}=\mathfrak{n}_{c}\oplus\mathfrak{n}_{n}$. The spin module $\mathcal{S}$ can be written as
 $$\mathcal{S}\cong\mathcal{S}_{\mathfrak{h}^{*}}\otimes\bigwedge\mathfrak{n}_{c}\otimes\bigwedge\mathfrak{n}^{*}_{n}\otimes\text{det}(\mathfrak{n}_{c}^{*})^{1/2}\otimes\text{det}(\mathfrak{n}_{n})^{1/2},$$
for some multiplicity space $\mathcal{S}_{\mathfrak{h}^{*}}$. It follows that $\mathcal{S}$ is a direct sum of copies of an irreducible $\mathfrak{g}$-representation of lowest weight $-\tilde{\rho}=\rho_{c}-\rho_{n}$ and its lowest weight space $S_{-\tilde{\rho}}$ is the tensor product of $\mathcal{S}_{\mathfrak{h}^{*}}$ and the lowest weight line in $\bigwedge\mathfrak{n}_{c}\otimes\bigwedge\mathfrak{n}^{*}_{n}\otimes\text{det}(\mathfrak{n}_{c}^{*})^{1/2}\otimes\text{det}(\mathfrak{n}_{n})^{1/2}$ ([FHT2] Lemma 1.15). Moreover on $L_{-\Lambda}\otimes S_{-\tilde{\rho}}\subset\mathcal{H}_{\Lambda}\otimes\mathcal{S}$ the following relations hold:
$$\gamma^{\alpha_{c}}=R_{-\alpha_{c}}=\sigma_{-\alpha_{c}}=0$$ for $\alpha_{c}\in\Delta_{c}^{+}$ and $$\gamma^{-\alpha_{n}}=R_{\alpha_{n}}=\sigma_{\alpha_{n}}=0$$ for $\alpha_{n}\in\Delta_{n}^{+}$. By direct calculation one finds (c.f. e.g. the proof of [FHT2] Proposition 1.19)
$$\restr{\big(\slashed D_{\mu}\big)_{c}+\big(\slashed D_{\mu}\big)_{n}}{L_{-\Lambda}\otimes S_{-\tilde{\rho}}}= i\gamma(\Lambda+\tilde{\rho}+\mu),$$ and so 
$$\restr{\big(\slashed D_{\mu}\big)_{c}^{2}+\big(\slashed D_{\mu}\big)_{n}^{2}}{L_{-\Lambda}\otimes S_{-\tilde{\rho}}}=|\Lambda+\tilde{\rho}+\mu|^2.$$ 
By another direct calculation, the operator $\big(\slashed D_{\mu}\big)_{c}^{2}+\big(\slashed D_{\mu}\big)_{n}^{2}$ commutes with the action of the Lie algebra and is therefore a scalar. Consequently 
\begin{align*}
 \slashed D_{\mu}^2 &=-\big(\slashed D_{0}\big)_{c}^{2}+\big(\slashed D_{0}\big)_{n}^{2}+ 2iT(\mu)-|\mu|^2\\
& =\big(\slashed D_{0}\big)_{c}^{2}+\big(\slashed D_{0}\big)_{n}^{2} -2 \big(\slashed D_{0}\big)_{c}^{2}+ 2iT(\mu)-|\mu|^2\\
&=|\Lambda+\tilde{\rho}|^2 -2\big(\slashed D_{0}\big)_{c}^{2}+ 2iT(\mu)-|\mu|^2\\
&=|\Lambda+\tilde{\rho}|^2 - \big(\slashed D_{0}\big)_{c}^{2}- \big(\slashed D_{\mu}\big)_{c}^{2}
\end{align*}
Note that $\big(\slashed D_{\mu}\big)^{2}_{c}$ is the Dirac operator for the compact group $K$. On each $K$-type $V_{K,\Lambda'}$ it is a nonnegative operator with a minimum eigenvalue of $|\Lambda'+\mu|^2$ achieved on the lowest weight line $L_{-\Lambda'}\subset V_{K,\Lambda'}$ ([FHT2] Proposition 1.19). Therefore
$$\restr{\slashed D_{\mu}^2}{V_{K,\Lambda'}\otimes\mathcal{S}}=|\Lambda+\tilde{\rho}|^2-|\Lambda'+\tilde{\rho}|^2- \big(\slashed D_{\mu}\big)^{2}_{c}.$$
But $-\Lambda$ being the minimal $K$-type means exactly that
$$|\Lambda+\tilde{\rho}|^2-|\Lambda'+\tilde{\rho}|^2\leq0,$$
which implies that $\slashed D_{\mu}^2$ is a non-negative operator whose minimum eigenvalue occurs on $L_{-\Lambda}\otimes S_{-\tilde{\rho}}$, the lowest weight space of the lowest $K$-type   in $\mathcal{H}\otimes\mathcal{S}$. The minimum eigenvalue is 
$$|\Lambda+\tilde{\rho}|^2-|\Lambda+\tilde{\rho}|^2- |\Lambda+\tilde{\rho}+\mu|^{2}$$
which is zero precisely when $\mu=-\Lambda-\tilde{\rho}$.
\end{proof}

\section{The Rossman Formula}

\textit{Throughout this section, $G$ is a connected real semisimple Lie group and ($\mathcal{H}_{-\Lambda},\pi)$ is an irreducible discrete series representation with lowest $K$-type $-\Lambda$.}

Let $U$ be an open neighborhood of $0\in\mathfrak{g}$ for which $\text{exp}: U\rightarrow G$ is invertible and analytic. For a compactly supported smooth function $\varphi\in C^{\infty}_{c}(U)$ the operator $\pi(\varphi):=\int_{\mathfrak{g}}dX\ \varphi(X)\pi(e^{X})$ is trace class and it is the aim of this section to establish the formula

$$\text{Tr}_{\mathcal{H}_{-\Lambda}}(\pi(\varphi))=\int_{\mathfrak{g}^{*}}\delta_{O_{-\Lambda-\tilde{\rho}}}(\xi)\ \Big(\int_{\mathfrak{g}} e^{i\langle\xi,X\rangle}\varphi(X)\hat{\mathcal{A}}(X)dX\Big)d\xi.$$
The proof is analogous the the compact case. 
One considers the trivial $\mathbb{Z}/2$-graded Hilbert bundle over $\mathcal{E}^{*}$ with fiber $\mathcal{H}_{\Lambda}\otimes\mathcal{S}$ and the formula of Definition \ref{family} defines a two-parameter family of superconnections on it. 
Only two points need to be addressed. The first is that $A_{\epsilon}^{t}(\varphi)$ and $\vec{\mu}(\varphi)$ are to be interpreted as the linear extensions $A_{\epsilon}^{t}(X)$ and $\vec{\mu}(X)$. 
The equivariant Chern characters are then defined in terms of $A_{\epsilon}^{t}(\varphi)$. The second point is that since the boundary of $\mathcal{E}^{*}$ is a little more interesting than the boundary of $\mathfrak{g}^{*}$, a little more care is required to establish invariance under one parameter deformations. The extra care comes in the form of the following lemma. 
Let $\mathcal{E}^{*}_{r}$ denote the intersection of $S(r):=\{|\xi|_{0}=r\}$\footnote{This is the length as measured in the positive definite metric gotten from reversing the sign of the Killing form $\langle\cdot,\cdot\rangle $ on its negative-definite subspace.} and $\mathcal{E}^{*}$. 
Let $\partial\mathcal{E}^{*}$ denote the boundary of $\mathcal{E}^{*}$ in $\mathfrak{g}^*$ and let  $\partial\mathcal{E}^{*}_r$ denote the intersection of $\partial\mathcal{E}^{*}$ with the interior of $S(r)$. Let $P(1)$ denote the line segment $(1,0)\to(\epsilon,0)$, $P(2)$ the line segment $(1,0)\to(1,1)$, and $P(3)$ the line segment $(1,1)\to(\epsilon,\epsilon)$. 

\begin{lem}\label{extracare}
Let $s$ be a parameter along each of the line segments above. Let $B_1$, $B_2$, and $B_3$ be the equivariant antiderivatives of $A_{\epsilon}^{t}(\varphi)$ along each of the line segments $P(1)$, $P(2)$, and $P(3)$ as given by Lemma \ref{deformation}. Then
$B_1=\text{str}(\frac{i}{2\sqrt{s}}\gamma(\xi)e^{A_{s}^{0}(\varphi)^2})$,
$B_2=\text{str}(\frac{1}{2\sqrt{s}}\slashed D_{0}e^{A^{s}_{1}(\varphi)^2})$, and
$B_3=\text{str}(\frac{1}{2\sqrt{s}}(\slashed D_{0}+i\gamma(\xi))e^{A_{s}^{s}(\varphi)^2})$, and all three can be extended-by-zero to $\partial\mathcal{E}^{*}$.
\end{lem}
\begin{proof}
The formulas for $B_{i}$ come directly from Lemma \ref{deformation}.
 When $\epsilon$ and $t$ are nonzero, as they are on all of $P(3)$ and a dense subset of $P(2)$, $B_i$ contains a Gaussian factor of the form $e^{-|\lambda+\xi|^2}$, peaked on the elliptic orbit of $\lambda$ and decaying in the stabilizer directions by equivariance. These stabilizer directions are orthogonal to the orbit in the Killing form and in these directions the boundary $\partial\mathcal{E}^{*}$ is infinitely far away from $\mathcal{E}^{*}$.\footnote{A good picture to keep in mind is the case of $SL_{2}(\mathbb{R})$, where $\mathfrak{g}^*$ can be identified with $\mathbb{R}^3$ such that $\mathcal{E}^{*}=\{(x,y,z)|z^2>x^2+y^2\}$, $\partial\mathcal{E}^{*}=\{(x,y,z)|z^2=x^2+y^2\}$. An elliptic orbit is a connected component of a two-sheeted hyperboloid opening in, say, the positive $z$-direction and asymptoting to $\partial\mathcal{E}^{*}$ and the stabilizer of a point on the hyperboloid is the line through that point and the origin, which only approaches the cone $\partial\mathcal{E}^{*}$ (which consists entirely of nilpotent elements) as the point moves to $\infty$ on the orbit.} Hence $B_2$ and $B_3$ extend-by-zero to $\partial\mathcal{E}^{*}$. Since $|\xi|^2=0$ on $\partial\mathcal{E}^{*}$ the formula for $B_1$ simplifies to  
 $$\frac{ie^{-|\gamma|^2}}{2\sqrt{s}}\text{str}(\gamma(\xi)e^{i\sqrt{s}d\gamma(\xi)+\vec{\mu}(\varphi)}),$$
 and since $d\gamma(\xi)=\gamma(\xi)d\xi$ any homogeneous term of positive degree (c.f. Remark \ref{topdeg}) will contain a factor of $\gamma(\xi)^2=-|\xi|^2=0$. So $B_1$ also extends-by-zero to $\partial\mathcal{E}^{*}$.
 \end{proof} 

To finish the proof we apply Lemma \ref{deformation} and Stokes theorem to the paths $P(1),P(2),$ and $P(3)$ to get
\begin{align*} 
&\int_{\mathcal{E}^{*}}\text{str}(e^{A_{\epsilon}^{0}(\varphi)^2}-e^{A_{1}^{0}(\varphi)^2})=\lim_{r\to\infty}\int_{\mathcal{E}^{*}_{r}\cup\partial\mathcal{E}^{*}_r}d\xi \int_{1}^{\epsilon}B_1ds,\\
&\int_{\mathcal{E}^{*}}\text{str}(e^{A_{1}^{1}(\varphi)^2}-e^{A_{1}^{0}(\varphi)^2})=\lim_{r\to\infty}\int_{\mathcal{E}^{*}_{r}\cup\partial\mathcal{E}^{*}_r}d\xi \int_{0}^{1}B_2,\\
&\int_{\mathcal{E}^{*}}\text{str}(e^{A_{\epsilon}^{\epsilon}(\varphi)^2}-e^{A_{1}^{1}(\varphi)^2})=\lim_{r\to\infty}\int_{\mathcal{E}^{*}_{r}\cup\partial\mathcal{E}^{*}_r}d\xi \int_{1}^{\epsilon}B_3ds.
\end{align*}
In all three cases, the integral over $\mathcal{E}^{*}_{r}$ vanishes in the $r\to\infty$ limit due to the Gaussian decay proved in Lemma \ref{decay}. Indeed, let $x_1,...,x_n$ be the standard coordinates on $\mathbb{R}^n$. Define $\langle \mathbf{x},\mathbf{y}\rangle=x_1y_1+...,x_ry_r-x_{r+1}y_{r+1}-...-x_ny_n$ and write $S(r)$ for the elements of length $r$ in the standard dot product. Then for each fixed $r$ the integrals mentioned above can be bounded by 
$$C\int_{\{\mathbf{x}\in\mathbb{R}^n|\langle\mathbf{x},\mathbf{x}\rangle>0\}\cap S(r)}e^{-\langle\mathbf{x},\mathbf{x}\rangle}d\mathbf{x}$$ where $C$ is a constant independent of $r$. The integral above is easily checked to vanish as $r\to\infty$.
In all three cases the integral over $\partial\mathcal{E}^{*}$ vanishes by Lemma \ref{extracare}. Similar to before, taking $\epsilon\to\infty$ in the segments $P(1)$ and $P(3)$ recovers the two sides of the Rossman formula and the proof is complete.

\section*{References}

\

[BGV]  N. Berline, E. Getlzer, M. Vergne: Heat Kernels and Dirac Operators \textit{Grundlehren Text Edition}, 2004

\

[FHT1]  D. Freed, M. Hopkins, and C. Teleman: Loop Groups and Twisted K-Theory \textit{Journal of Topology} \textbf{4}, 2011

\

[FHT2]  D. Freed, M. Hopkins, and C. Teleman: Loop Groups and Twisted K-Theory \textit{Journal of the American Mathematical Society} \textbf{26}, 2013

\

[FHT3]  D. Freed, M. Hopkins, and C. Teleman: Loop Groups and Twisted K-Theory \textit{Annals of Math} \textbf{174}, 2011.

\

[FT] D. Freed and C. Teleman: Dirac Families for Loop Groups as Matrix factorizations \textit{Comptes Rendus de l'Acadademie des Sciences} Paris, \textbf{353}, (2015)

\

[Ki]  N. Kitchloo: Dominant $K$-theory and integrable highest-weight representations of Kac-Moody groups \textit{Advances in Mathematics} \textbf{221}, 2009

\

[Kn]  A. Knapp: Representation Theory of Semisimple Groups \textit{Princeton Landmarks in Mathematics}, 1997

\

[N]  Y. Nakagama: Spectral Analysis in Krein Spaces \textit{Publications of the Research Institute for Mathematical Sciences, Kyoto University} \textbf{24}, 1988 

\

Evans Hall

University Drive, Berkeley, CA 94720

kiranluecke@berkeley.edu

\end{document}